\documentclass[12pt]{article}
\usepackage{cite}
\usepackage{amsmath, amsthm, amssymb,slashed}
\usepackage{ifpdf}
\ifpdf
  \usepackage[pdftex]{graphicx}
  \usepackage{epstopdf}
\else
  \usepackage[dvips]{graphicx}
\fi
\parskip=\baselineskip
\def\Bbb{\mathbb}
\def\Tr{{\rm Tr}}
\def\16{{\bf 16}}
\def\1{{\bf 1}}
\def\2{{\bf 2}}
\def\3{{\bf 3}}
\def\4{{\bf 4}}
\def\t{{\mathbf t}}

\def\U{{\mathcal U}}
\def\Bbb{\mathbb}
\def\frak{\mathfrak}
\def\C{{\Bbb {C}}}
\def\Z{{\Bbb Z}}
\def\A{{\mathcal A}}
\def\CS{{\mathrm{CS}}}
\def\d{{\mathrm d}}
\def\U{{\mathcal U}}
\def\Tr{{\mathrm{Tr}}}
\def\R{{\Bbb{R}}}
\def\D{{\mathcal D}}

\def\R{{\Bbb{R}}}\def\Z{{\Bbb{Z}}}
\def\N{{\mathcal N}}

\font\teneurm=eurm10 \font\seveneurm=eurm7 \font\fiveeurm=eurm5
\newfam\eurmfam
\textfont\eurmfam=\teneurm \scriptfont\eurmfam=\seveneurm
\scriptscriptfont\eurmfam=\fiveeurm

\font\teneusm=eusm10 \font\seveneusm=eusm7 \font\fiveeusm=eusm5
\newfam\eusmfam
\textfont\eusmfam=\teneusm \scriptfont\eusmfam=\seveneusm
\scriptscriptfont\eusmfam=\fiveeusm

\font\tencmmib=cmmib10 \skewchar\tencmmib='177
\font\sevencmmib=cmmib7 \skewchar\sevencmmib='177
\font\fivecmmib=cmmib5 \skewchar\fivecmmib='177
\newfam\cmmibfam
\textfont\cmmibfam=\tencmmib \scriptfont\cmmibfam=\sevencmmib
\scriptscriptfont\cmmibfam=\fivecmmib

\numberwithin{equation}{section}

\def\d{\mathrm d}
\def\C{{\Bbb C}}
\def\Z{{\Bbb Z}}
\def\A{{\mathcal A}}

\begin{document}
\begin{titlepage}
\begin{flushright}
hep-th/yymm.nnnn
\end{flushright}
\vskip 1.5in
\begin{center}
{\bf\Large{Khovanov Homology And Gauge Theory}}
\vskip
0.5cm {Edward Witten} \vskip 0.05in {\small{ \textit{School of
Natural Sciences, Institute for Advanced Study}\vskip -.4cm
{\textit{Einstein Drive, Princeton, NJ 08540 USA}}}
}
\end{center}
\vskip 0.5in
\baselineskip 16pt
\abstract{In these notes, I will sketch a new approach to Khovanov homology of knots and links
based on counting the solutions of certain elliptic partial differential equations in four and five dimensions.
The equations are formulated on four and five-dimensional
manifolds with boundary, with a rather subtle boundary condition that encodes the knots and links.  The construction is formally analogous to Floer and Donaldson theory in three and four dimensions.  It was 
discovered using quantum field theory
arguments but can be described and understood purely in terms of classical gauge theory.}
\date{September, 2011}
\end{titlepage}
\def\Hom{\mathrm{Hom}}

\section{Introduction}

The goal of these lecture notes is to propose how certain elliptic differential equations in four
and five dimensions can be used
for a new description of the Jones polynomial and Khovanov homology.
My inspiration was a previous physics-based approach to
Khovanov homology \cite{GSV}.  I sought to find a gauge theory alternative to that construction.

In the presentation here, we  will achieve brevity  by omitting all the ``why'' 
questions. In other words, we will explain what construction in four- and five-dimensional gauge theory is
supposed to reproduce the Jones polynomial and Khovanov homology, but not how this
answer emerges from dualities of quantum fields and strings.  For such questions,
along with much more detail and extensive references to related work, we refer to
\cite{fiveknots}.  That reference also contains much more detail concerning the ``what'' questions;
our explanations here are rather schematic.

Without relying on quantum field theory arguments, can one
directly deduce a known construction of the Jones polynomial or Khovanov homology from
the gauge theory framework that will be described here?  A strategy to do this
for the Jones polynomial has been
outlined in \cite{GWnew}, and some of the ideas  will be summarized 
at the end of these notes.  There is not yet an analogous computation for Khovanov homology.

\section{Gradient Flow Of Complex Chern-Simons}
We start by remembering the relation of the Chern-Simons function in three dimensions
to the instanton equation in four dimensions.  In this talk, 
$G$ is a compact simple Lie group and $A$ is a connection on 
a $G$-bundle $E\to W$, where $W$ is some manifold.  We write $G_\C$ for 
the 
complexification of $G$, and $\A$ for a connection on a $G_\C$ bundle, such as
 the complexification $E_\C$ of $E$. We write $\U$ and $\U_C$ for the spaces of, respectively,
 connections on $E$ or on $E_\C$.  Finally, by an elliptic equation, we mean an equation 
 that is elliptic modulo the action of the gauge group.

 If $W$ is a three-manifold, then a connection $A$ on the $G$-bundle $E\to W$
has a Chern-Simons invariant
\begin{equation}\CS(A)=\frac{1}{4\pi}\int_W\Tr\,\left(A\wedge \d A +\frac{2}{3}A\wedge A\wedge A\right).\end{equation}
To arrive at the instanton equations, we pick a Riemannian metric on $W$ and then place
the obvious Riemannian metric on the space $\U$ of connections:
\begin{equation}\d s^2=-\int_W\Tr\,\delta A\wedge \star \delta A.\end{equation}
Then, viewing $-\CS(A)$ as a Morse function on $\U$, we write the equation of gradient flow:
\begin{equation}
\frac{\d A}{\d s}=\nabla  \CS(A).
\end{equation}
Something nice happens; the equation of gradient flow turns out to have four-dimensional symmetry.  It is equivalent to the instanton equation on the four-manifold
$M=W\times \R$:
\begin{equation}F^+=0.\end{equation}  (For a two-form $\omega$ in four dimensions,
we write $\omega^+$ and $\omega^-$ for its selfdual and anti-selfdual projections.) This fact is the starting point for Floer cohomology of three-manifolds and its relation to Donaldson theory of four-manifolds.

We want to do the same thing, roughly speaking, for the complex Lie group $G_\C$.
To begin with, a connection $\A$ on a $G_\C$ bundle $E_\C\to W$ has a Chern-Simons
function: 
\begin{equation}\CS(\A)=\frac{1}{4\pi}\int_W\Tr\,\left(\A\wedge\d \A+
\frac{2}{3}\A\wedge\A\wedge\A\right).\end{equation}
 To do Morse theory, we have to make two immediate changes.  
 First, a Morse function is supposed to be real, but $\CS(\A)$ is 
 actually complex-valued.  So we pick a complex number  $e^{i\alpha}$ 
 of modulus 1 and define our Morse function to be (provisionally) 
\begin{equation}h_0=-\mathrm{Re}\,(e^{i\alpha}\CS(\A)).\end{equation}

Second, there is not any convenient metric on the space  $\U_\C$ of complex connections that
has the full $G_\C$ gauge symmetry.  So we pick a Kahler metric on 
$\U_\C$ that is invariant only under $G$, not $G_\C$:
\begin{equation}\d s^2=-\int_W\Tr\,\delta \A\wedge \star \delta\overline \A.\end{equation}
Now we can write a gradient flow equation
\begin{equation}\frac{\d \A}{\d s}=-\nabla h_0.\end{equation}

However, we are really usually interested in a complex 
connection $\A$ up to complex-valued gauge transformations, 
but here we have written an equation that is only invariant under
unitary gauge transformations.   To compensate for this, we should set the moment map
to zero and consider the previous equation only in the space of zeroes of the moment map.  
In other words, if we decompose $\A$ in real and imaginary parts as 
$\A=A+i\phi$, where $A$ is a real connection and 
$\phi\in \Omega^1(W)\otimes \mathrm{ad}(E)$, then the Kahler manifold $\U_\C$ has a Kahler form
\begin{equation}\omega=\int_W\Tr\,\delta A\wedge \star \delta \phi.\end{equation}
The moment map for the action of $G$-valued gauge transformations is
\begin{equation} \mu=\d_A\star\phi,\end{equation}
and we should really consider the previous gradient flow equations 
in the space of zeroes of the moment map.

However, it is somewhat better to introduce another field $\phi_0$
as a sort of Lagrange muliplier.  $\phi_0$ is a section of the real adjoint bundle $\mathrm{ad}(E)$
and we write an extended Morse function 
\begin{equation}h=h_0+\int_W\d^3x \sqrt g \,\Tr\,\phi_0 \mu\end{equation}
whose critical points are all at $\mu=0$.  On the space of $\phi_0$ fields, we place the
obvious metric
\begin{equation}\d s^2=-\int_W\d^3x\sqrt g \,\Tr \,\delta \phi_0^2,\end{equation}
 and now, writing $\Phi$ for the pair $(\A,\phi_0)$, we write the gradient flow equations
\begin{equation}\frac{\d\Phi}{\d s}=-\nabla h(\Phi).\end{equation}

Something nice happens, just like what happened in the real case. 
The flow equations in this sense are elliptic partial differential equations with a full four-dimensional symmetry.   They can be written
\begin{align}\notag   (F-\phi\wedge\phi)^+ &= t(\d_A\phi)^+\\
\notag (F-\phi\wedge\phi)^-&=-t^{-1}(\d_A\phi)^-\\
\d_A\star\phi & = 0,\label{zelf}\end{align}
with
\begin{equation}t=\frac{1-\cos\alpha}{\sin\alpha}.\end{equation}
These are elliptic differential equations modulo the action of the gauge group, 
for each $t\in \Bbb{RP}^1=\R\cup\infty$ (for $t\to 0$ or $\infty$, multiply the second 
equation by $t$ or the first by $t^{-1}$).

In writing the equations, we have combined the imaginary part of the connection, 
$\phi\in \Omega^1(W)\otimes \mathrm{ad}(E)$, with the Lagrange multiplier $\phi_0\in \mathrm{ad}(E)$, to a field (also called $\phi$) that takes values in 
$\Omega^1(M)\otimes\mathrm{ad}(E)$.  We are using the fact that for $M=W\times \R$, we have
$T^*M=T^*W\oplus \R$, where $\R$ is a trivial real line bundle of rank 1.

The equations (\ref{zelf}) can be the starting point for developing a Floer-like theory for
the complex Lie group $G_\C$.  That is the case we will want for the Jones polynomial
and Khovanov homology.  However, before focusing on that case, let us ask what might one
mean by Floer theory for other real forms of $G$.

If $t=0$ or $\infty$, the equations (\ref{zelf}) admit the involution $\phi\to -\phi$.  If we simply
ask for a solution to be invariant under this operation, we will set $\phi$ to zero and reduce
to the selfdual or anti-selfdual Yang-Mills equations (depending on whether $t$ is 0 or $\infty$).
This will lead to ordinary Floer theory.  To get something new, we can instead pick an involution $\tau$
(ether inner or outer) of the compact Lie group $G$.  Conjugacy classes of such involutions are in natural correspondence with real forms of $G$. Imposing invariance under $\phi\to -\phi$
together with the involution $\tau$, we get a reduced set of equations appropriate to constructing
a Floer-like theory for the real form that corresponds to $\tau$.  In this way,
one can hope to construct a Floer-like theory for any real form of $G$.

Let us return to the full system of equations (\ref{zelf}), without reducing to a real form.
These equations were actually first studied in the context of gauge theory and geometric
Langlands \cite{KW}.   Roughly speaking, one constructs a family of 
four-dimensional topological field theories
that are just like Donaldson theory except based on the equations 
(\ref{zelf}) instead of the instanton equations.  Geometric Langlands duality is then naturally formulated as an equivalence between the theories that arise at two different values of $t$.
From that vantage point, geometric Langlands duality is a consequence of $S$-duality of 
$\N=4$ super Yang-Mills theory in four dimensions.  (I am leaving out numerous details, one of
which is that it is natural to consider complex values of $t$.)

It does not seem likely that this theory gives interesting 
four-manifold invariants, and there may be a technical problem in defining them, 
analogous to what happens in Donaldson theory on a four-manifold of $b_2^+=1$.  
In geometric Langlands duality, four-manifold invariants are not the point.
One is mainly interested not in invariants of four-manifolds but in the structures 
that this same topological field theory attaches to two-manifolds (categories 
of boundary conditions) and
to three-manifolds (spaces of  physical states or ``morphisms'').  These are 
essentially not affected by the technical difficulties that may affect the four-manifold invariants.

Although this theory probably does not give interesting four-manifold invariants, 
the literature on geometric Langlands  does give a reason to believe in hindsight, 
and perhaps even in foresight,
that it gives interesting knot invariants.  An extension of the usual geometric Langlands duality
that is known as ``quantum'' geometric Langlands is related to quantum groups
\cite{gaits}, which in turn can be used to describe the Jones polynomial and similar invariants of knots.  
Possibly this should have made us anticipate that the topological field theory associated to 
geometric Langlands can be used to compute the Jones polynomial.  Another possible
clue relating Khovanov homology to geometric Langlands can be gleaned from 
an attempt to describe Khovanov homology as an $A$-model \cite{ss} and more specifically
from the relation of the symplectic varieties used in that attempt to spaces of Hecke modifications
\cite{ck}.  The framework that we sketch below indeed leads to spaces of Hecke modifications, though of
a somewhat different type.

\section{Relation To The Jones Polynomial}

To make contact with the Jones polynomial, we take the four-manifold $M$ to be $W\times \R_+$, where $W$ is
a three-manifold and $\R_+$ is a half-line $y\geq 0$.  For $y\to\infty$, we require the fields to approach a chosen
critical point (or more precisely a critical gauge orbit) of the Morse function.  The equations (\ref{zelf}) were chosen so that a critical point is a complex flat connection,
corresponding to a homomorphism $\rho:\pi_1(W)\to G_\C$.  ($\rho$ must satisfy a mild condition
of semistability \cite{Corlette} so that the moment map condition $\d_A\star\phi=0$ can be obeyed.)  The simplest case
is that $W$ is $\sf S^3$ or $\R^3$ so that there is only one possible critical point, with $A=\phi=0$ up to
a $G$-valued gauge transformation.

\begin{figure}
 \begin{center}
   \includegraphics[width=3.5in]{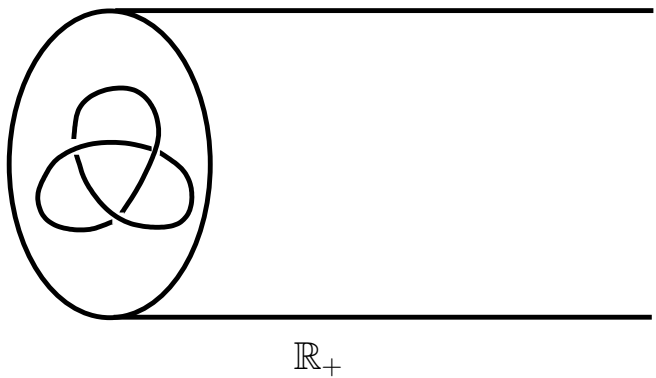}
 \end{center}
\caption{\small  A knot has been placed at the boundary of the four-manifold $M=W\times\R_+$. }
 \label{knotboundary}
\end{figure}
Meanwhile, as sketched in fig. \ref{knotboundary}, knots are placed at $y=0$, that is at the endpoint of $\R_+$.  
The boundary condition at $y=0$ is an elliptic boundary condition that is a little involved to explain, and we will 
postpone it to section \ref{bc}.  But this boundary condition depends on the choice of a knot (or link\footnote{In the case
of a link $L=\cup_i K_i$ with components $K_i$, each component can be labeled by a different representation
$R_i$ of $G^\vee$.}) $K$, and on the choice of a representation
$R$ of the Langlands or GNO dual group $G^\vee$ of $G$.  In this description, this is the only way that $K$
(or $R$) enters. 

For the next step, we imitate Donaldson theory.  Let $a_n$ be the ``number'' of solutions of the equations (\ref{zelf})
with instanton number $n$.  (The count of solutions is made in a standard sort of way; for example, an isolated
nondegenerate solution is weighted by a sign, which is the sign of the fermion determinant.)  Conjecturally,
$a_n$ vanishes for large enough $|n|$, but this has not been proved.  
 The instanton number is
defined in the usual way as a multiple of 
\begin{equation}\int_M\Tr\,F\wedge F,\end{equation}
where $F=\d A+A\wedge A$ is the curvature.  Actually, to define the instanton number as a topological invariant
is subtle, because $M=W\times \R_+$ has ``ends'' at $y=0,\infty$.   The boundary conditions at the two ends
give 
identifications of the bundle $E$ with fixed bundles (namely a specified flat bundle at $y=\infty$ and an embedding of the
tangent bundle of $W$ at $y=0$; see section \ref{bc} for the latter), and this enables one to define the instanton 
number as a topological invariant once a framing of $W$ is picked.  The  dependence on the 
framing ultimately reproduces the framing dependence of the Jones polynomial and its generalizations.

Then we introduce a variable $q$ and define the counting function
\begin{equation}\label{gold} J(q;K,R)=\sum_n a_n q^n.\end{equation}
For $G^\vee=SU(2)$ and $R$ the two-dimensional representation, and with $W=\R^3$, 
the claim is that this function is the Jones polynomial.
In general, for $W=\R^3$, we expect to get the usual knot and link invariants associated to quantum groups and Chern-Simons gauge theory.

I refer to \cite{fiveknots,newlook} for  explanations of why $J(q;K,R)$ is expected to reproduce the Jones polynomial
and its generalizations.   Here I will just mention 
the main steps.  First, one represents the function $J(q;K,R)$ by a path integral in $\N=4$ super Yang-Mills theory.
This can be done in a standard way because the equations (\ref{zelf}) whose solutions we want 
to count are associated to one of the
supersymmetry generators of the $\N=4$ theory.  The second step is electric-magnetic duality -- 
the same quantum symmetry of $\N=4$ super Yang-Mills theory that underlies geometric Langlands duality.  
This puts the path integral
in a more convenient form.  The third and most novel step is to interpret the path integral of the $\N=4$ theory on
$W\times\R_+$ as an unusual sort of three-dimensional path integral.  The most familiar quantum field 
theory representation of the Jones polynomial \cite{wittenjones} is as an integral over the 
space $\U$ of real-valued connections:\footnote{The object $\Tr_R\,P\exp\oint_KA$ is the trace of the holonomy
of the connection $A$ around the loop $K$, with the trace taken in the representation $R$.}
\begin{equation}\label{zoz}\int_\U DA\,\exp\left(ik \CS(A)\right)\,\Tr_R\,P\exp\oint_KA. \end{equation}
We can generalize this to an integral over a middle-dimensional cycle $\Gamma$ in $\U_\C$, the space of
complex-valued connections:
\begin{equation}\label{oz}\int_\Gamma D\A\,\exp\left(ik\CS(\A)\right)\,\Tr_RP\exp\oint_K\A.  \end{equation}
This is a generalization of (\ref{zoz}), to which it reduces if $\Gamma=\U$.  In \cite{fiveknots,newlook}, it is shown
that the $\N=4$ path integral on $W\times\R_+$ (in its electric-magnetic dual version), is equivalent to a 
generalized Chern-Simons path integral (\ref{oz}) with $\Gamma$ the Lefschetz thimble associated to  the critical
point that was used to define the boundary condition at $y=\infty$.  In establishing this relation, a key fact
is that the equations (\ref{zelf}) whose solutions we want to count are the gradient flow equations of complex
Chern-Simons theory.  For a generic $W$, the Lefschetz thimble $\Gamma$ is not equivalent to the real
integration cycle $\U$ of the usual Chern-Simons path integral (\ref{zoz}), so the counting function 
(\ref{gold}) cannot be directly compared to the invariants derived from quantum groups and 
Chern-Simons theory.  However, precisely for $W=\R^3$,
they are equivalent, and that is why in this case, the counting of solutions is expected to reproduce the 
Jones polynomial and its generalizations. 

\section{Khovanov Homology}

However, our goal is Khovanov homology, not the Jones polynomial, 
This means we are supposed to ``categorify'' the situation, and associate to a knot or link 
a vector space rather than a number.    A suitable trace in the vector space 
will give back the original numerical invariant, which in the present context is the counting function of eqn. (\ref{gold}).
In plain words, this means that the approach to the Jones polynomial
just described has to be derived
from a picture in one more dimension.  Physicists would describe the situation as follows: viewing the extra 
dimension as ``time,'' quantization gives a physical ``Hilbert'' space, which will be the Khovanov homology, and 
then if we compactify the extra dimension on a circle, we get a trace leading back to the original theory -- in this
case, the theory that computes the Jones polynomial. 

Let us practice by ``categorifying'' something -- in fact, we will categorify the Casson invariant.  The 
Casson invariant is an invariant of a three-manifold $W$.   It is defined by 
``counting'' (up to gauge transformation and with suitables signs) the
flat connections $A$ on a $G$-bundle $E\to W$. A flat connection 
is a solution of the equation
\begin{equation}F=0.\end{equation}
In three dimensions, this is not an elliptic equation, but is part of a nonlinear elliptic 
complex. Just as for linear elliptic complexes,  it is sometimes convenient to 
``fold'' the complex and reduce to the case of an ordinary elliptic equation, 
rather than a complex.   In the present example, this is done by
introducing a section $\phi_0$ of $\mathrm{ad}(E)$ and replacing the equation 
$F=0$ with the Bogomolny equation
\begin{equation}\label{bog}F+\star \d_A\phi_0=0.\end{equation}
The count of solutions is the same, apart from rather degenerate solutions, since a simple vanishing theorem 
says that any smooth solution on a compact manifold $W$ has 
$\phi_0$ covariantly constant.  But now the equation is elliptic.  This gives a more convenient starting point for
categorification.

``Categorification'' of the Casson invariant is accomplished simply by   replacing 
$\phi_0$ by the covariant
derivative with respect to a new coordinate, $\phi_0\to D/Ds$.  Thus we replace the three-manifold
$W$ by the four-manifold $X=\R\times W$, where $\R$ is parametrized by the ``time'' $s$, and
we substitute $\phi_0\to D/Ds$.   This makes sense, in that the substitution gives 
a differential equation on $X$ (rather than a differential operator), because 
we started with an equation in which
$\phi_0$ only appears inside the commutator $\d_A\phi_0=[\d_A,\phi_0]$. 
This commutator is now
replaced by $[\d_A,D/Ds]$, which is a component of the four-dimensional curvature.  

Generically, a procedure like this, even if it gives a differential equation, 
will not give an elliptic one, let alone one with four-dimensional symmetry.  In this case, 
however, we actually get back the instanton equation $F^+=0$, in a slightly different way from the way we got it 
before.  What follows from 
this is that the Casson invariant -- a numerical invariant  
computed by counting solutions of the original equation $F=0$ -- can be categorified to 
Floer cohomology.

Let us spell this out in detail.  To define Floer cohomology, one introduces a chain complex
that has a basis corresponding to the flat connections on $W$.  The differential in the chain
complex is determined by counting suitable instanton solutions on $X=\R\times W$.  However,
we can compute the Euler characteristic of the Floer cohomology without having to know the differential:
we simply sum over the  basis vectors of the chain complex, weighted by appropriate signs.  
  This sum
gives back the Casson invariant.  So the Casson invariant is the Euler characteristic of the Floer cohomology.
The fancy way to say this is that Floer cohomology is a categorification of the Casson invariant.

So far there is no novelty here; this construction is essentially due to Floer.  Now we want to similarly 
categorify the Jones polynomial, which from the point of 
view of the present lecture is the invariant associated to counting solutions of the equations  (\ref{zelf}) with certain boundary conditions.  On a generic four-manifold $M$, we would 
have no way to proceed as there is no candidate for a field $\phi_0$ that will be replaced by 
$D/Ds$. 

However, if $M=W\times \R_+$, which is the case 
if we are studying the Jones polynomial in the way  suggested here, then 
$T^*M=T^*W\oplus \R$, where the part of $\phi$ associated
to the second summand is the field $\phi_0$ that we originally introduced as a Lagrange multiplier. 
We categorify by introducing a 
new dimension and replacing $\phi_0\to D/Ds$, as before. 

In this way, we get a partial differential equation on the five-manifold
$X=\R\times W\times \R_+$.  Moreover, this turns out to be an elliptic equation. (See \cite{fiveknots} for details.  The
same equation was also introduced in \cite{haydys}.)
And if we set $t=1$, we get a full four-dimensional symmetry; that is, the five-dimensional
equation  can be naturally formulated on
$X=M^*\times \R_+$ for any oriented four-manifold $M^*$.

\begin{figure}
 \begin{center}
   \includegraphics[width=3.5in]{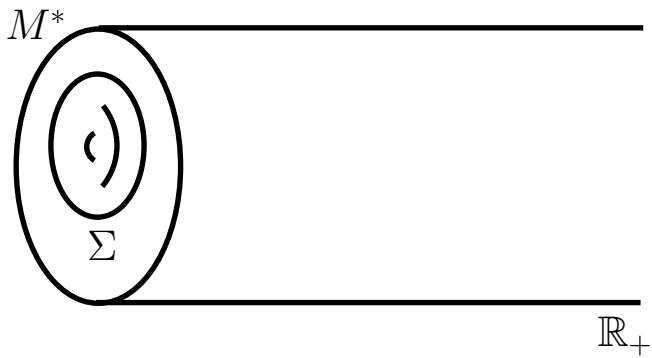}
 \end{center}
\caption{\small  A two-manifold $\Sigma$ embedded in the boundary of $X=M^*\times \R_+$. }
 \label{secondknot}
\end{figure}
The boundary condition  at $y=0$ (described in section \ref{bc}) that can be used to get the Jones polynomial    can be ``lifted'' to five dimensions, roughly by 
$\phi_0\to D/Ds$.  The boundary is now a four-manifold $M^*$ rather than a three-manifold $W$.   
Instead of modifying the boundary condition along a knot $K\subset W$, we now
modify it along a two-manifold $\Sigma\subset M^*$, as in fig. \ref{secondknot}.

To get the candidate for Khovanov homology,  we specialize to the time-independent case
$M^*=\R\times W$, $\Sigma=\R\times K$.   Then, following Floer (and following, for example, \cite{KM}, where
a similar strategy is used to categorify the Alexander polynomial), we define a chain complex which has a 
basis given by the time-independent solutions, that is the solutions of the four-dimensional equations
\begin{align}\notag   (F-\phi\wedge\phi)^+ &= t(\d_A\phi)^+\\
\notag (F-\phi\wedge\phi)^-&=-t^{-1}(\d_A\phi)^-\\
\d_A\star\phi & = 0.\end{align}  
The differential in the chain complex is constructed in a standard fashion by counting certain time-dependent 
solutions.  (Here we use the fact that the five-dimensional equations can themselves be interpreted in terms of 
gradient flow.)  The cohomology of this differential is the candidate for Khovanov homology. 

The candidate Khovanov homology is $\Z\times \Z$-graded, like the real thing, where one
grading is the cohomological grading, and the second grading, sometimes called the $q$-grading,  
is associated to the instanton number, integrated
over $W\times \R_+$.    (Because  $W\times \R_+$ is not compact and has a boundary, the 
definition of the instanton number has subtleties that turn out to match the framing anomalies of 
Chern-Simons theory.  See section \ref{bc}.)  One can define an equivariant Euler characteristic of the candidate Khovanov homology.
In passing to the Euler characteristic, one loses the cohomological grading, but one keeps the grading by
instanton number.  So the equivariant Euler characteristic is a function of a single variable $q$.  Because the complex that computes the candidate Khovanov homology has a basis given
by the same time-independent solutions that have to be counted to compute the function $J(q;K,R)$ defined
in eqn. (\ref{gold}), the equivariant
Euler characteristic is equal to this function.  This corresponds to the usual relation between Khovanov homology
and the Jones polynomial. The reasoning involved is exactly parallel to the relation between Floer cohomology
and the Casson invariant.

\begin{figure}
 \begin{center}
   \includegraphics[width=2in]{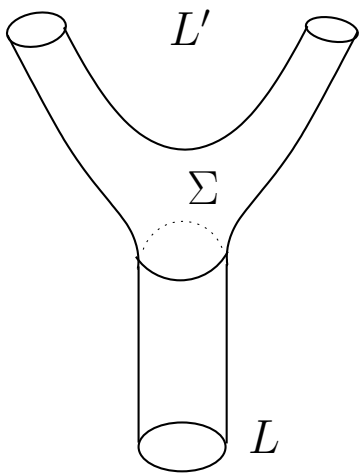}
 \end{center}
\caption{\small  The surface $\Sigma$ represents a ``link cobordism'' from an initial link $L$ (in this case, an unknot)
to a final link $L'$ (in this case, two disjoint unknots).  We embed $\Sigma$ in the boundary of $X=M^*\times \R_+$, as in fig. \ref{secondknot}.}
 \label{cob}
\end{figure}

We are not limited to the case that the five-manifold $X$ and the Riemann surface $\Sigma$  are time-independent.  Preserving the time-independence of $X$ but allowing
a more general $\Sigma$, as in fig. \ref{cob},  we get candidates for the ``link cobordisms'' of Khovanov homology.
Explicitly, to such a picture, we associate a linear transformation from the Khovanov homology of the initial link to that of the final
link.  It is  computed by counting five-dimensional solutions with time-dependent
boundary conditions (determined by $\Sigma$), asymptotic in the past and future to specified
time-independent solutions that represent initial and final states.

\section{The Boundary Condition}\label{bc}

In this section, we will describe the boundary condition that is needed for this construction, initially in the absence of knots. 
 It is essentially enough to describe the boundary condition in four dimensions rather than five 
 (once one understands it, the lift to five dimensions is fairly obvious),   and as the boundary condition 
 is local, we assume initially that the boundary of the four-manifold is just $\R^3$.
So we work on $M=\R^3\times \R_+$.   (This special case is anyway the right case for
the Jones polynomial, which concerns knots in $\R^3$ or equivalently $\sf S^3$.)

We will have to make use of one of the important equations in gauge theory, which is Nahm's equation.  
Nahm's equation is a system of ordinary differential equations for a triple $X_1,X_2,X_3$ valued in $\frak g^3$, 
where $\frak g$ is the Lie algebra of $G$.
The equations read
\begin{equation}\frac{\d X_1}{\d y}+[X_2,X_3]=0\end{equation}
and cyclic permutations. On a half-line $y\geq 0$, Nahm's equations have the special solution
$$ X_i=\frac{\mathbf t_i}{y},$$
where the $\mathbf t_i$ are elements of $\frak g$ that obey the $\frak{su}(2)$  commutation relations
$[\t_1,\t_2]=\t_3$, etc.  We are mainly interested in the case that the $\t_i$ define a ``principal
$\frak{su}_2$ subalgebra'' of $\frak g$, in the sense of Kostant.

This sort of singular solution of Nahm's equations was important
in the work of Nahm on monopoles, and in later work of Kronheimer and others.   We will
use it to define an elliptic boundary condition for our equations.   In fact, 
Nahm's equations can be embedded in our four-dimensional
equations on $\R^3\times \R_+$.   If we look for a solution that is (i) invariant under translations of
$\R^3$, (ii) has the connection $A=0$, (iii) has $\phi=\sum_{i=1}^3\phi_i\,\d x_i + 0\cdot \d y$
(where $x_1,x_2,x_3$ are coordinates on $\R^3$ and $y$ is the normal coordinate), then our four-dimensional equations 
(\ref{zelf}) reduce to Nahm's equations
\begin{equation}\frac{\d \phi_1}{\d y}+[\phi_2,\phi_3]=0,\end{equation}
and cyclic permutations.  So the ``Nahm pole'' gives a special solution of our equations
\begin{equation}\phi_i=\frac{\t_i}{y}.\end{equation}
We  define an elliptic boundary condition by declaring that we will allow only solutions
that are asymptotic to this one for $y\to 0$.

This is the boundary condition that we want at $y=0$, in the absence of knots, and for $W=\R^3$. For the most obvious boundary condition for getting Khovanov homology, we require that $A,\phi\to 0$ for $y\to\infty$. 
It is plausible (but unproved) that with these conditions, the special solution with 
the Nahm pole is the only one.  (This would correspond to Khovanov 
homology of the empty knot being of rank 1.)

Something interesting happens when we generalize this boundary condition from $W=\R^3$ to an arbitrary oriented
Riemannian three-manifold $W$.  We describe the situation first for $G=SO(3)$.  To make sense for general $W$ of the formula $\phi_i=\t_i/y$
or $\phi=\sum_i \t_i \,\d x^i/y$ (where the $x^i$ are local coordinates on $W$), we have to interpret the
objects $\t_i$ as supplying an isomorphism between the bundle $\mathrm{ad}(E)$ and the tangent bundle of $W$.
A study of the equations (\ref{zelf}) shows further that this isomorphism must be covariantly constant, so the gauge
field $A$, when restricted to $W=\partial M$, can be identified with the Riemannian connnection of $W$.  It follows
from this together with some elementary arguments that to define the instanton number, integrated over $W\times \R_+$,
as a topological invariant, one needs  a framing of $W$ (a trivialization of its tangent bundle).\footnote{One also needs
a condition at $y=\infty$, that is at the infinite end of $\R_+$.  This is provided by the requirement that the
connection approaches a specified critical point for  $y\to \infty$.}  For $G=SU(2)$, the
tangent bundle of $W$ has to be lifted to an $SU(2)$ bundle, so 
a choice of spin structure on $W$ is part of the boundary condition.  The analog of this for any compact $G$ is
that, when restricted to $W$, the connection $A$ is obtained from the Riemannian connection on $W$ via
a principal embedding $\rho:\frak{su}(2)\to \frak g$.  A choice of spin structure is part of the boundary condition 
precisely if
the subgroup of $G$ associated to such a principal embedding is $SU(2)$ rather than $SO(3)$.

The next step is to modify the boundary condition to incorporate the presence of a knot $K$.
 The local model is that the boundary is $\R^3$, and $K$ is a copy of $\R\subset \R^3$.
 The boundary condition is described by giving a 
 singular model solution on $\R^3\times \R_+$ that along the boundary 
 has the now-familiar Nahm pole away from $K$, but has some other behavior 
 along $K$.  There should be one model solution for every irreducible representation $R$ of the dual
 group $G^\vee$.  The model solution is invariant under translations along $K$, 
 so it can be obtained by solving some reduced equations on $\R^2\times \R_+$.  Before explaining a little
 more about this, we will first explain another reason that it is important to study the reduced equations in three
 dimensions.

 \section{Categories}\label{categ}
 
 \begin{figure}
 \begin{center}
   \includegraphics[width=3.5in]{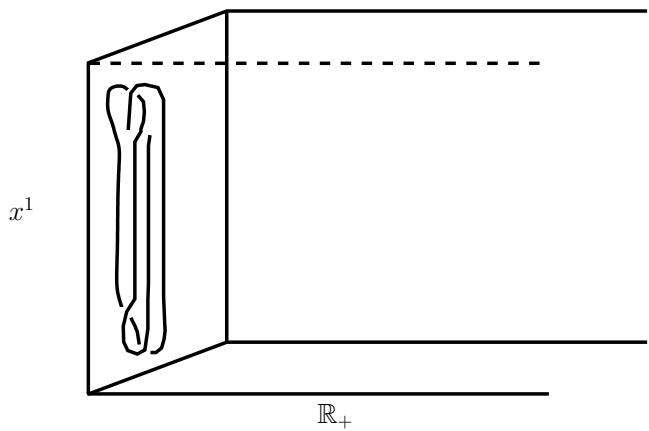}
 \end{center}
\caption{\small Stretching a knot in one dimension, to reduce to a description in one dimension less.  One of the
directions in $\R^3$ -- here labeled as $x^1$ -- plays the role of ``time.''  After stretching of the knot,
one hopes that a solution of the equations becomes almost everywhere nearly independent of $x^1$. If so,
a knowledge of the solutions that are actually independent of $x^1$ can be a starting point for understanding
four-dimensional solutions.  This type of analysis will fail at the critical points of the function $x^1$ along a knot or link;
a correction has to be made at those points.}
 \label{stretch}
\end{figure}
   To compute the Jones polynomial,  
 we need to count certain
 solutions in four dimensions; knowledge of these solutions is also the first step 
 in constructing the candidate for Khovanov homology.   How are we supposed
 to describe four-dimensional solutions?  A standard strategy, often used in Floer theory
 and its cousins, involves ``stretching'' the knot in one direction (fig. \ref{stretch}), in the hope of reducing to
 a piecewise description by solutions in one dimension less.  To carry out such a strategy, we need to be
 able to understand the reduced equations in three dimensions.

  \begin{figure}
 \begin{center}
   \includegraphics[width=3.5in]{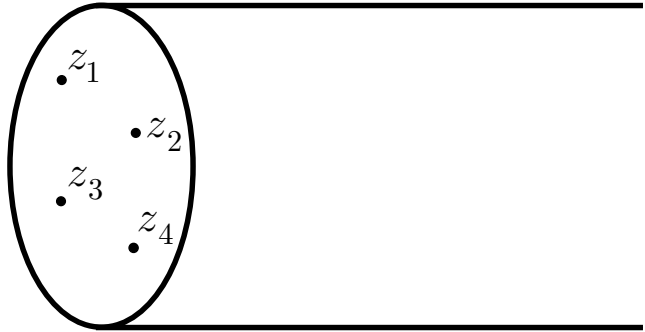}
 \end{center}
\caption{\small To construct Khovanov homology, one wishes to attach a suitable category to $\R^2$ or $\sf S^2$
with marked points labeled by representations of $G^\vee$.  In the framework described here, to construct
this category one has to consider the reduced equations on $\R^2\times \R_+$ (or $\sf S^2\times \R_+$, as sketched
here), with marked points on the boundary.}
 \label{category}
\end{figure}
 Another way to make the point is as follows.   Most 
  mathematical definitions of Khovanov homology proceed, directly or implicitly, by 
  defining a category of objects associated to a two-sphere (or in some versions, a 
  copy of $\C=\R^2$) with marked points that are suitably labeled.   
  In the present approach, this category should be the $A$-model category of the moduli 
  space of solutions of the reduced three-dimensional equations in the appropriate geometry,
  as sketched in fig. \ref{category}.  (The $A$-model in question is really an $A$-model with a superpotential;
  it is described in \cite{GWnew} as a Landau-Ginzburg $A$-model.)
 
 The equations when reduced to three dimensions
 have a pleasantly simple 
 structure.    After getting  
 a simplification via a small vanishing theorem 
 for some of the fields, the equations can be schematically described as follows. There are three operators $\D_i$ (constructed from $A$ and $\phi$) that commute,
 \begin{equation} [\D_i,\D_j]=0,~~i,j=1,2,3.\end{equation}
 And they obey a ``moment map'' constraint
 \begin{equation}\sum_{i=1}^3[\D_i,\D_i{}^\dagger]=0.\end{equation}
 
 The construction of the $\D_i$ in terms of $A$ and $\phi$ depends on $t$.  At $t=1$,
 the equations  just described coincide with the ``extended Bogomoly equations'' \cite{KW},
 which describe geometric Hecke transformations in the context of the geometric Langlands correspondence.  
 For this value of $t$, the requisite model solutions, with singularities labeled by representations of the dual group,
  can be described in closed form \cite{fiveknots,mikhaylov}.
The appearance of geometric Hecke transformations is not a complete surprise, since one of the known descriptions of Khovanov homology \cite{ck} is based on a 
$B$-model category of moduli spaces of
 geometric Hecke transformations.    (This approach has not yet been derived from quantum field theory,
 though some of the steps are clear.)

 For generic $t$, a more transparent structure actually emerges \cite{GWnew}, even though, in this
 case, closed form solutions are not known.  
 After a change of variables that is only valid for $t\not=\pm 1$,  the equations 
 $[\D_i,\D_j]=0=\sum_i[\D_i,\D_i^\dagger]$  are actually more familiar.  They describe a flat
 $G_\C$ bundle $E\to\R^2\times\R_+$ endowed with a hermitian metric that obeys a moment map
 condition.  For a special value of $t$, the moment map condition is the one studied
 long ago in \cite{Corlette}.  As far as is known, the $t$-dependence of the moment map is not important
 for any qualitative results.   
 
 Since $\R^2\times\R_+$ is simply-connected, how can we get
 anything interesting from a flat connection?  The answer is that there is additional
 structure in the behavior at $y=0$ (and $\infty$).  A flat bundle over $\R^2\times \R_+$ is, of course, the pullback
of a flat bundle on $\R^2$, which we will think of as $\C$.   The boundary conditions
at $y=0$ gives the flat bundle $E\to \C$ the structure of an ``oper,'' in the language of
geometric Langlands.   At the points corresponding to the knots, the oper has
singularities, but the flat bundle has no monodromy around these singularities. Such oper 
singularities are classified again by representations of the dual group, as reviewed in \cite{frenkel}.   
This gives a way to understand for generic $t$ the existence of model solutions with singularities labeled
by a representation of the dual group.  (For generic $t$, these solutions are not known in closed form, but
they can be found numerically for $G^\vee=SU(2)$.)
As also explained in \cite{frenkel}, opers with
singularities of this type are related to integrable models such as the Gaudin spin chain.  

Integrable spin systems are part of the context in which the Jones polynomial was originally discovered.  
So the relation of the reduced three-dimensional equations to an integrable spin chain offers some encouragement
that we may be able to directly understand how the counting of the solutions of the four-dimensional equations
(\ref{zelf}) leads to the Jones polynomial.  Actually carrying this out \cite{GWnew} requires one more
key idea, which we describe next. 

\section{Knot Projections}
\begin{figure}
 \begin{center}
   \includegraphics[width=2in]{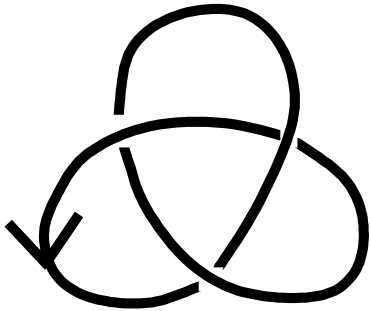}
 \end{center}
\caption{\small  Projection of a knot to two dimensions.}
 \label{projection}
\end{figure}

In standard approaches to the Jones polynomial and Khovanov homology, one often begins by considering a
projection of a knot to two dimensions, as in fig. \ref{projection}.  One makes some definition based on this projection,
and then one endeavors to prove that this definition is actually independent of the choice of a projection. 

In the context of the four- and five-dimensional equations considered here, there is a very natural way to incorporate
a  knot projection by modifying the boundary condition for $y\to\infty$.  For brevity, we will describe this
only for four-dimensional gauge theory on $M=\R^3\times \R_+$.  Instead of requiring that $A,\phi\to 0$ for $y\to\infty$,
we keep that condition on $A$, but we change the condition on $\phi$.
We pick a triple $c_1,c_2,c_3$ of commuting
elements of $\frak t$, the Lie algebra of a maximal torus $T\subset G$,  and we ask for
\begin{equation}\phi\to\sum_i c_i\cdot \d x^i\end{equation}
for $y\to\infty$. Here $x^1,x^2,x^3$ are Euclidean coordinates on $\R^3$.
We use the fact that the equations have an exact solution for $A=0$
and $\phi$ of the form just indicated.   There is a nice theory of Nahm's equations with this sort of 
condition at infinity \cite{Kron}.
 
 The counting of solutions of an elliptic equation is constant
under continuous variations (provided certain conditions are obeyed) so one expects
that the Jones polynomial can be computed with this more general asymptotic
condition, for an arbitrary choice of $\vec c=(c_1,c_2,c_3)$.  

If $G=SU(2)$, then $\frak t$ is one-dimensional.  
So if $\vec c$ is non-zero, it has the form $\vec c=c\cdot\vec a$ where $c$ is
a fixed (nonzero) element of $\frak t$ and $\vec a$ is a vector in three-space. 
So picking
$\vec c$ essentially means picking 
a vector $\vec a$ pointing in some direction in $\R^3$. 
The choice of $\vec a$ determines a projection of $\R^3$ to a plane, so this
is now built into the construction.
 For $G$ of higher
rank, one could do something more general, but it seems sufficient to take
$\vec c=c \cdot \vec a$ with $c$ a regular element of $\frak t$.    The advantage to this method of introducing a knot
projection is that we know {\it a priori} that what we compute will be independent of the choice of the projection.  

Taking $\vec c\not=0$ is described by physicists as
``gauge symmetry 
breaking'' or ``moving on the Coulomb branch.''  A closely parallel construction is important in the 
theory of superconductivity, 
in Taubes's proof that ``SW=GW,'' and in the theory of weak interactions and Higgs particles.
 (For a review of these matters, emphasizing their common features, see \cite{BAMS}.)     Taking $\vec c$ sufficiently
generic
gives a drastic simplification because the equations become quasi-abelian in a  certain sense.
 On a length scale large than $1/|\vec c|$, the solutions can be almost 
everywhere approximated by solutions of an abelian version of the same equations.   
There is an important locus where this fails, but it can be understood.

By reasoning along these lines, it is possible \cite{GWnew} to give a fairly direct map from  the counting 
function $J(q;K,R)$, defined in eqn. (\ref{gold}),  to a known construction of
the Jones polynomial (and its analogs for other representations and
gauge groups). The known construction that has been obtained this way is the vertex model (which is presented, for example,
in \cite{kauffman}). One of the main tools in the analysis was the relationship of the reduced three-dimensional
equations to integrable models, as indicated in section \ref{categ}.   The analysis gave a way of relating the
counting function $J(q;K,R)$ to the Jones polynomial, without relying on arguments of quantum field theory.

\vskip 2 cm \noindent
{\bf{Acknowledgment:}} Research supported in part by NSF Grant PHY-0969448. 
\bibliographystyle{unsrt}

\end{document}